\documentclass[12pt,lenq]{amsart}
\usepackage{amsmath}
\usepackage{amscd}
\usepackage{amssymb}
\usepackage{amsbsy}
\usepackage{amsfonts}
\usepackage{latexsym}
\usepackage{graphics}
\usepackage{amsmath,amscd,latexsym}
\usepackage{multirow}
\usepackage{array}
\usepackage{paralist}
\usepackage{titletoc}

\theoremstyle{plain}
\newtheorem{theorem}{Theorem}[section]
\newtheorem{proposition}[theorem]{Proposition}

\newtheorem{corollary}[theorem]{Corollary}

\newtheorem{definition}[theorem]{Definition}

\newtheorem{main theorem}[theorem]{Main Theorem}

\newtheorem{question}[theorem]{Question}

\newtheorem{conjecture}[theorem]{Conjecture}

\newlength\savewidth

\newcommand{\ZZ}{\mathbb{Z}}
\newcommand{\QQ}{\mathbb{Q}}
\newcommand{\RR}{\mathbb{R}}
\newcommand{\CC}{\mathbb{C}}
\newcommand{\HH}{\mathbb{H}}
\newcommand{\QQQ}{\hat{\mathbb{Q}}}
\newcommand{\RRR}{\hat{\mathbb{R}}}

\newcommand{\Conway}{\mbox{\boldmath$S$}^{2}}
\newcommand{\Conways}
{(\mbox{\boldmath$S$}^{2},\mbox{\boldmath$P$})}
\newcommand{\PP}{\mbox{\boldmath$P$}}
\newcommand{\PConway}{\mbox{\boldmath$S$}}
\newcommand{\ptorus}{\mbox{\boldmath$T$}}
\newcommand{\OO}{\mbox{\boldmath$O$}}

\newcommand{\PSL}{\mbox{$\mathrm{PSL}$}}
\newcommand{\SL}{\mbox{$\mathrm{SL}$}}
\newcommand{\tr}{\mbox{$\mathrm{tr}$}}
\newcommand{\Isom}{\mbox{$\mathrm{Isom}$}}

\newcommand{\DD}{\mathcal{D}}
\newcommand{\RGPP}[1]{\hat\Gamma_{#1}}
\newcommand{\RGP}[1]{\Gamma_{#1}}
\newcommand{\curve}{\mathcal{C}}
\newcommand{\plamination}{\mathcal{PL}}
\newcommand{\einv}{\mathcal{E}}

\newcommand{\svert}{\,|\,}

\newcommand{\llangle}{\langle\langle}
\newcommand{\rrangle}{\rangle\rangle}

\newcommand{\lp}{(\hskip -0.07cm (}
\newcommand{\rp}{)\hskip -0.07cm )}

\begin{document}

\title{Simple loops on 2-bridge spheres in 2-bridge link complements}

\author{Donghi Lee}
\address{Department of Mathematics\\
Pusan National University \\
San-30 Jangjeon-Dong, Geumjung-Gu, Pusan, 609-735, Republic of Korea}
\email{donghi@pusan.ac.kr}

\author{Makoto Sakuma}
\address{Department of Mathematics\\
Graduate School of Science\\
Hiroshima University\\
Higashi-Hiroshima, 739-8526, Japan}
\email{sakuma@math.sci.hiroshima-u.ac.jp}

\subjclass[2000]{Primary 57M25, 20F06 \\
\indent {The first author was supported by
a 2-Year Research Grant of Pusan National University.
The second author was supported
by JSPS Grants-in-Aid 22340013 and 21654011.}}


\begin{abstract}
The purpose of this note is to announce complete answers to
the following questions.
(1) For an essential simple loop
on a 2-bridge sphere in a 2-bridge link complement,
when is it null-homotopic in the link complement?
(2) For two distinct essential simple loops
on a 2-bridge sphere in a 2-bridge link complement,
when are they homotopic in the link complement?
We also announce applications of these results
to character varieties and
McShane's identity.
\end{abstract}
\maketitle


\section{Introduction}
Let $K$ be a knot or a link in $S^3$
and $S$ a punctured
sphere in the complement $S^3-K$
obtained from a bridge sphere of $K$.
Then the following natural question arises.

\begin{question} \label{question0}
{\rm (1)} Which essential simple loops on $S$ are null-homotopic in
$S^3-K$?

{\rm (2)} For two distinct essential simple loops on $S$,
when are they
homotopic in $S^3-K$?
\end{question}

A refined version of the first question for
$2$-bridge spheres of $2$-bridge links
was proposed
in the second author's joint work with
Ohtsuki and Riley \cite[Question ~9.1(2)]{Ohtsuki-Riley-Sakuma},
in relation with epimorphisms
between $2$-bridge links.
It may be regarded as a special variation of a question
raised by Minsky ~\cite[Question ~5.4]{Gordon}
on essential simple loops on Heegaard surfaces of $3$-manifolds.

The purpose of this note is to announce a complete answer
to Question ~\ref{question0} for
$2$-bridge spheres of $2$-bridge links
established by the series of papers
\cite{lee_sakuma, lee_sakuma_2, lee_sakuma_3, lee_sakuma_4}
and to explain its application
to the study of character varieties and McShane's identity
\cite{lee_sakuma_5}.

The key tool for solving the question is small cancellation theory,
applied to two-generator and one-relator presentations
of $2$-bridge link groups.
We note that it has been proved
by Weinbaum ~\cite{Weinbaum} and
Appel and Schupp ~\cite{Appel-Schupp}
that the word and conjugacy problems
for prime alternating link groups are solvable,
by using small cancellation theory
(see also \cite{Johnsgard} and references in it).
Moreover, it was shown
by Sela ~\cite{Sela} and Pr\'eaux ~\cite{Preaux} that
the word and conjugacy problems
for any link group are solvable.
A characteristic feature of our work
is that it gives complete answers to
special (but also natural)
word and conjugacy problems
for the link groups of 2-bridge links,
which form a special (but also important) family
of prime alternating links.
(See \cite{Adams, ASWY} for the role of $2$-bridge links
in Kleinian group theory.)

This note is organized as follows.
In Sections ~\ref{statements}, \ref{application1} and \ref{application2},
we describe the main results,
applications to character varieties and McShane's identity.
The remaining sections are devoted to explanation of
the idea of the proof of the main results.
In Section ~\ref{group_presentation},
we describe the two-generator and one-relator
presentation of the $2$-bridge link group
to which small cancellation theory is applied,
and give a natural decomposition of the relator,
which plays a key role in the proof.
In Section ~\ref{sequences},
we introduce a certain finite sequence
associated with the relator and state its key properties.
In Section ~\ref{small_cancellation},
we recall small cancellation theory and
present a characterization of the ``pieces''
of the symmetrized subset arising from the relator.
In Sections ~\ref{outline_proof} and \ref{outline_proof2},
we describe outlines of the proofs of
the main results.

The authors would like to thank
Norbert A'Campo,
Hirotaka Akiyoshi,
Brian Bowditch, Danny Calegari,
Max Forester,
Koji Fujiwara,
Yair Minsky,
Toshihiro Nakanishi,
Caroline Series
and Ser Peow Tan for stimulating conversations.

\section{Main results}
\label{statements}

\begin{figure}[h]
\begin{center}
\includegraphics{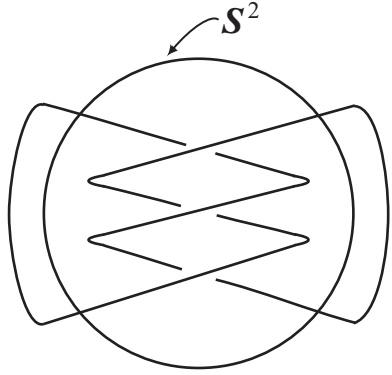}
\end{center}
\caption{
\label{fig.trefoil}
$(S^3,K(r))=(B^3,t(\infty))\cup (B^3,t(r))$ with $r=1/3$.
Here $(B^3,t(r))$ and $(B^3,t(\infty))$, respectively,
are the inside and the outside of the bridge sphere $\Conway$.
}
\end{figure}

For a rational number $r \in \QQQ:=\QQ\cup\{\infty\}$,
let $K(r)$ be the $2$-bridge link of slope $r$,
which is defined as the sum
$(S^3,K(r))=(B^3,t(\infty))\cup (B^3,t(r))$
of rational tangles of slope $\infty$ and $r$
(see Figure ~\ref{fig.trefoil}).
The common boundary $\partial (B^3,t(\infty))= \partial (B^3,t(r))$
of the rational tangles is identified
with the {\it Conway sphere} $\Conways:=(\RR^2,\ZZ^2)/H$,
where $H$ is the group of isometries
of the Euclidean plane $\RR^2$
generated by the $\pi$-rotations around
the points in the lattice $\ZZ^2$.
Let $\PConway$ be the $4$-punctured sphere $\Conway-\PP$
in the link complement $S^3-K(r)$.
Any essential simple loop in $\PConway$, up to isotopy,
is obtained as
the image of a line of slope $s\in\QQQ$ in
$\RR^2-\ZZ^2$
by the covering projection onto $\PConway$.
The (unoriented) essential simple loop in $\PConway$ so obtained
is denoted by $\alpha_s$.
We also denote by $\alpha_s$ the conjugacy class of
an element of $\pi_1(\PConway)$
represented by (a suitably oriented) $\alpha_s$.
Then the {\it link group} $G(K(r)):=\pi_1(S^3-K(r))$
is identified with
$\pi_1(\PConway)/ \llangle\alpha_{\infty},\alpha_r\rrangle$,
where $\llangle\cdot \rrangle$ denotes the normal closure.

\begin{figure}[h]
\begin{center}
\includegraphics{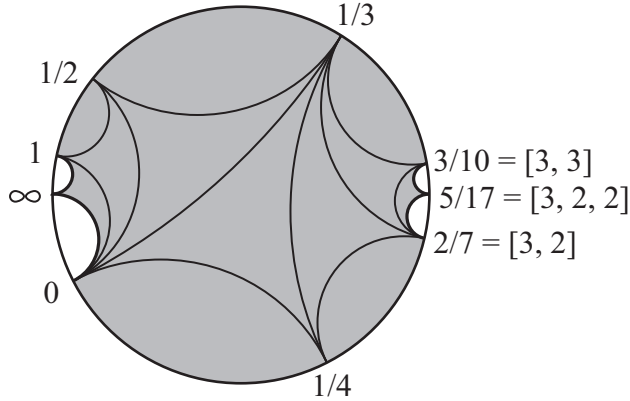}
\end{center}
\caption{\label{fig.fd}
A fundamental domain of $\hat\Gamma_r$ in the
Farey tessellation (the shaded domain) for $r=5/17=[3,2,2]$.}
\end{figure}

Let $\DD$ be the
{\it Farey tessellation},
whose ideal vertex set is identified with $\QQQ$.
For each $r\in \QQQ$,
let $\RGP{r}$ be the group of automorphisms of
$\DD$ generated by reflections in the edges of $\DD$
with an endpoint $r$, and
let $\RGPP{r}$ be the group generated by $\RGP{r}$ and $\RGP{\infty}$.
Then the region, $R$, bounded by a pair of
Farey edges with an endpoint $\infty$
and a pair of Farey edges with an endpoint $r$
forms a fundamental domain of the action of $\RGPP{r}$ on $\HH^2$
(see Figure ~\ref{fig.fd}).
Let $I_1$ and $I_2$ be the closed intervals in $\RRR$
obtained as the intersection with $\RRR$ of the closure of $R$.
Suppose that $r$ is a rational number with $0<r<1$.
(We may always assume this except when we treat the
trivial knot and the trivial $2$-component link.)
Write
\begin{center}
\begin{picture}(230,70)
\put(0,48){$\displaystyle{
r=
\cfrac{1}{m_1+
\cfrac{1}{ \raisebox{-5pt}[0pt][0pt]{$m_2 \, + \, $}
\raisebox{-10pt}[0pt][0pt]{$\, \ddots \ $}
\raisebox{-12pt}[0pt][0pt]{$+ \, \cfrac{1}{m_k}$}
}} \
=:[m_1,m_2, \dots,m_k],}$}
\end{picture}
\end{center}
where $k \ge 1$, $(m_1, \dots, m_k) \in (\mathbb{Z}_+)^k$, and $m_k \ge 2$.
Then the above intervals are given by
$I_1=[0,r_1]$ and $I_2=[r_2,1]$,
where
\begin{align*}
r_1 &=
\begin{cases}
[m_1, m_2, \dots, m_{k-1}] & \mbox{if $k$ is odd,}\\
[m_1, m_2, \dots, m_{k-1}, m_k-1] & \mbox{if $k$ is even,}
\end{cases}\\
r_2 &=
\begin{cases}
[m_1, m_2, \dots, m_{k-1}, m_k-1] & \mbox{if $k$ is odd,}\\
[m_1, m_2, \dots, m_{k-1}] & \mbox{if $k$ is even.}
\end{cases}
\end{align*}

We recall the following fact
(\cite[Proposition ~4.6 and Corollary ~4.7]{Ohtsuki-Riley-Sakuma}
and \cite[Lemma ~7.1]{lee_sakuma})
which describes
the role of $\RGPP{r}$ in the study of
$2$-bridge link groups.

\begin{proposition}
\label{previous_results}
{\rm (1)}
If two elements $s$ and $s'$ of $\QQQ$ belong to the same orbit $\RGPP{r}$-orbit,
then the unoriented loops $\alpha_s$ and $\alpha_{s'}$ are homotopic in $S^3-K(r)$.

{\rm (2)}
For any $s\in\QQQ$,
there is a unique rational number
$s_0\in I_1\cup I_2\cup \{\infty, r\}$
such that $s$ is contained in the  $\RGPP{r}$-orbit of $s_0$.
In particular, $\alpha_s$ is homotopic to $\alpha_{s_0}$ in
$S^3-K(r)$.
Thus if $s_0\in\{\infty, r\}$, then $\alpha_s$ is null-homotopic
in $S^3-K(r)$.
\end{proposition}

Thus the following question naturally arises
(see \cite[Question ~9.1(2)]{Ohtsuki-Riley-Sakuma}).

\begin{question}
\label{question}
{\rm (1)}
Does the converse to Proposition ~\ref{previous_results}(2) hold?
Namely, is it true that $\alpha_s$ is null-homotopic in $S^3-K(r)$
if and only if $s$ belongs to the $\RGPP{r}$-orbit of $\infty$ or $r$?

{\rm (2)}
For two distinct rational numbers $s, s' \in I_1\cup I_2$,
when are the unoriented loops $\alpha_s$ and $\alpha_{s'}$
homotopic in $S^3-K(r)$?
\end{question}

The following theorem proved in \cite{lee_sakuma}
gives a complete answer to Question ~\ref{question}(1).

\begin{theorem}
\label{first_result}
The loop $\alpha_s$ is null-homotopic in $S^3 - K(r)$
if and only if $s$ belongs to the $\RGPP{r}$-orbit of $\infty$ or $r$.
In other words, if $s\in I_1\cup I_2$, then
$\alpha_s$ is not null-homotopic in $S^3-K(r)$.
\end{theorem}

This has the following application to the study of
epimorphisms between $2$-bridge link groups
(see \cite[Section ~2]{lee_sakuma} for precise meaning).

\begin{corollary}
\label{first_cor}
There is an upper-meridian-pair-preserving epimorphism
from $G(K(s))$ to $G(K(r))$
if and only if $s$ or $s+1$
belongs to the $\RGPP{r}$-orbit of $r$ or $\infty$.
\end{corollary}

The following theorem
proved in \cite{lee_sakuma_2, lee_sakuma_3, lee_sakuma_4}
gives a complete answer
to Question ~\ref{question}(2).

\begin{theorem}
\label{second_result}
Suppose that $r$ is a rational number such that $0 < r \le 1/2$.
For distinct $s, s' \in I_1\cup I_2$,
the unoriented loops $\alpha_s$ and $\alpha_{s'}$ are
homotopic in $S^3-K(r)$
if and only if one of the following holds.
\begin{enumerate}[\indent \rm (1)]
\item
$r=1/p$, where $p \ge 2$ is an integer,
and $s=q_1/p_1$ and $s'=q_2/p_2$ satisfy
$q_1=q_2$ and $q_1/(p_1+p_2)=1/p$, where $(p_i, q_i)$ is a pair of
relatively prime positive integers.
\item
$r=3/8$, namely $K(r)$ is the Whitehead link,
and the set $\{s, s'\}$ equals
either $\{1/6, 3/10\}$ or $\{3/4, 5/12\}$.
\end{enumerate}
\end{theorem}

The proof of Theorem ~\ref{second_result}
reveals the structure of the normalizer
of an element of $G(K(r))$ represented by $\alpha_s$.
This enables us to show the following.

\begin{theorem}
\label{peripheral_results}
Let $r$ be a rational number such that $0 < r \le 1/2$.
Suppose $K(r)$ is hyperbolic,
i.e., $r=q/p$ and $q\not\equiv \pm1 \pmod{p}$,
and let $s$ be a rational number contained in $I_1\cup I_2$.

{\rm (1)} The loop $\alpha_s$
is peripheral
if and only if one of the following holds.
\begin{enumerate}[\indent \rm (i)]
\item
$r=2/5$ and $s=1/5$ or
$s=3/5$.
\item
$r=n/(2n+1)$ for some integer $n\ge 3$,
and $s=(n+1)/(2n+1)$.
\item
$r=2/(2n+1)$ for some integers $n \ge 3$,
and $s=1/(2n+1)$.
\end{enumerate}

{\rm (2)} The conjugacy class $\alpha_s$ is primitive in $G(K(r))$
with the following exceptions.
\begin{enumerate}[\indent \rm (i)]
\item
$r=2/5$ and $s=2/7$ or $3/4$.
In this case $\alpha_s$ is the third power of some primitive element
in $G(K(r))$.
\item
$r=3/7$ and $s=2/7$.
In this case $\alpha_s$ is the second power of some primitive element
in $G(K(r))$.
\item
$r=2/7$ and $s=3/7$.
In this case $\alpha_s$ is the second power of some primitive element
in $G(K(r))$.
\end{enumerate}
\end{theorem}

At the end of this section, we describe a relation of
Theorem ~\ref{first_result} with
the question raised by Minsky in \cite[Question ~5.4]{Gordon}.
Let $M=H_+\cup_S H_-$ be a Heegaard splitting of
a $3$-manifold $M$.
Let $\Gamma_{\pm}:=MCG(H_{\pm})$ be the mapping
class group of $H_{\pm}$, and
let $\Gamma^0_{\pm}$ be the kernel of the map
$MCG(H_{\pm})\to \mathrm{Out}(\pi_1(H_{\pm}))$.
Identify $\Gamma^0_{\pm}$ with a subgroup of $MCG(S)$,
and consider the subgroup $\langle \Gamma^0_+,\Gamma^0_-\rangle$
of $MCG(S)$.
Now let $\Delta_{\pm}$ be the set of (isotopy classes of)
simple loops in $S$ which bound a disk in $H_{\pm}$.
Let $Z$ be the set of essential simple loops in $S$
which are null-homotopic in $M$.
Note that $Z$ contains $\Delta_{\pm}$ and invariant
under $\langle \Gamma^0_+,\Gamma^0_-\rangle$.
In particular, the orbit
$\langle \Gamma^0_+,\Gamma^0_-\rangle(\Delta_+\cup\Delta_-)$
is a subset of $Z$.
Then Minsky posed the following question.

\begin{question}
\label{Minsky_question}
When is $Z$ equal to the orbit
$\langle \Gamma^0_+,\Gamma^0_-\rangle(\Delta_+\cup\Delta_-)$?
\end{question}

The above question makes sense
not only for Heegaard splittings but also
bridge decompositions of knots and links.
In particular, for $2$-bridge links, the groups
$\RGP{\infty}$ and $\RGP{r}$
in our setting correspond to the groups
$\Gamma^0_+$ and $\Gamma^0_-$, and hence the group
$\RGPP{r}$ corresponds to the group
$\langle\Gamma^0_+,\Gamma^0_-\rangle$.
To make this precise, recall the bridge decomposition
$(S^3,K(r))=(B^3,t(\infty))\cup (B^3,t(r))$,
and let $\tilde\Gamma_{+}$ (resp. $\tilde\Gamma_{-}$) be the mapping
class group of the pair $(B^3,t(\infty))$ (resp. $(B^3,t(r))$), and
let $\tilde\Gamma^0_{\pm}$ be the kernel of the natural map
$\tilde\Gamma_{+}\to \mathrm{Out}(\pi_1(B^3-t(\infty)))$
(resp. $\tilde\Gamma_{-}\to \mathrm{Out}(\pi_1(B^3-t(r)))$).
Identify $\tilde\Gamma^0_{\pm}$ with a subgroup of
the mapping class group $MCG(\PConway)$ of the
$4$-times punctured sphere $\PConway$.
Recall that the Farey tessellation $\DD$ is identified
with the curve complex of $\PConway$ and there is a natural epimorphism
from $MCG(\PConway)$ to the automorphism group $Aut(\DD)$
of $\DD$,
whose kernel is isomorphic to $(\ZZ/2\ZZ)^2$.
Then the groups $\Gamma_{\infty}$ and $\Gamma_r$, respectively,
are identified with the images of
$\tilde\Gamma^0_{+}$ and $\tilde\Gamma^0_{-}$
by this epimorphism.
Moreover, the sets $\{\alpha_{\infty}\}$
and $\{\alpha_{r}\}$, respectively,
correspond to the sets $\Delta_+$ and $\Delta_-$.
Theorem ~\ref{first_result} says that
the set $Z$ of simple loops in $\PConway$
which are null-homotopic in $S^3-K(r)$
is equal to the orbit
$\langle \Gamma_{\infty}, \Gamma_r\rangle(\Delta_+\cup\Delta_-)$.
Thus Theorem ~\ref{first_result} may
be regarded as an answer to
the special variation of Question ~\ref{Minsky_question}.

\section{Application to character varieties}
\label{application1}

In this section and the next section,
we assume $r=q/p$, where $p$ and $q$ are relatively prime
positive integers such that
$q\not\equiv \pm1 \pmod{p}$.
This is equivalent to the condition that
$K(r)$ is hyperbolic,
namely the link complement $S^3-K(r)$ admits
a complete hyperbolic structure of finite volume.
Let $\rho_r$ be the $\PSL(2,\CC)$-representation of $\pi_1(\PConway)$
obtained as the composition
\[
\pi_1(\PConway) \to
\pi_1(\PConway)/ \llangle\alpha_{\infty},\alpha_r\rrangle
\cong
\pi_1(S^3-K(r))
\to
\Isom^+(\HH^3)
\cong
\PSL(2,\CC),
\]
where the last homomorphism is the holonomy representation
associated with the complete hyperbolic structure.

Now, let $\ptorus$ be the once-punctured torus
obtained as the quotient
$(\RR^2-\ZZ^2)/\ZZ^2$,
and let $\OO$ be the orbifold
$(\RR^2-\ZZ^2)/\hat H$
where $\hat H$ is the group generated by $\pi$-rotations
around the points in $(\frac{1}{2}\ZZ)^2$.
Note that $\OO$ is the orbifold with
underlying space a once-punctured sphere
and with three cone points of cone angle $\pi$.
The surfaces $\ptorus$ and $\PConway$, respectively,
are $\ZZ/2\ZZ$-covering and $(\ZZ/2\ZZ)^2$-covering of $\OO$,
and hence their fundamental groups are identified
with subgroups
of the orbifold fundamental group $\pi_1(\OO)$
of indices $2$ and $4$, respectively.
The $\PSL(2,\CC)$-representation $\rho_r$ of $\pi_1(\PConway)$
extends, in a unique way, to that of $\pi_1(\OO)$
(see \cite[Proposition ~2.2]{ASWY}),
and so we obtain, in a unique way, a
$\PSL(2,\CC)$-representation of $\pi_1(\ptorus)$
by restriction.
We continue to denote it by $\rho_r$.
Note that $\rho_r:\pi_1(\ptorus) \to \PSL(2,\CC)$
is {\it type-preserving}, i.e.,
it satisfies the following conditions.
\begin{enumerate}
\item
$\rho_r$ is irreducible, i.e.,
its image does not have a common fixed point on
$\partial \HH^3$.
\item
$\rho_r$ maps a peripheral element of $\pi_1(\ptorus)$
to a parabolic transformation.
\end{enumerate}
By extending the concept of a geometrically infinite end
of a Kleinian group,
Bowditch ~\cite{Bowditch2}
introduced the notion of the end invariants of a type-preserving
$\PSL(2,\CC)$-representation of $\pi_1(\ptorus)$.
Tan, Wong and Zhang ~\cite{Tan_Wong_Zhang_6} (cf. \cite{Tan_Wong_Zhang_1})
extended this notion (with slight modification)
to an arbitrary $\PSL(2,\CC)$-representation of $\pi_1(\ptorus)$.
(To be precise, \cite{Tan_Wong_Zhang_6}
treats $\SL(2,\CC)$-representations.
However, the arguments work for $\PSL(2,\CC)$-representations.)

To recall the notion of end invariants,
let $\curve$ be the set of free homotopy classes
of essential simple loops on $\ptorus$.
Then $\curve$ is identified with $\QQQ$, the vertex set of
the Farey tessellation $\DD$ by the following rule.
For each $s\in\QQQ$,
let $\beta_s$ be the essential simple loop
on $\ptorus$ obtained as the image of a line
of slope $s$ in $\RR^2-\ZZ^2$.
Then the correspondence $s\mapsto \beta_s$
gives the desired identification $\QQQ\cong \curve$.
The projective lamination space $\plamination$
is then identified with $\RRR:=\RR\cup\{\infty\}$
and contains $\curve$ as the dense subset
of rational points.

\begin{definition}
\label{def_end_invariant}
{\rm
Let $\rho$ be a $\PSL(2,\CC)$-representation of $\pi_1(\ptorus)$.

(1) An element $X\in \plamination$ is an {\it end invariant}
of $\rho$
if there exists a sequence of distinct elements
$X_n\in\curve$ such that
$X_n\to X$ and such that $\{|\tr \rho(X_n)|\}_n$
is bounded from above.

(2) $\einv(\rho)$ denotes the set of end invariants of $\rho$.
}
\end{definition}

In the above definition, it should be noted that
$|\tr \rho(X_n)|$ is well-defined
though $\tr \rho(X_n)$ is defined only up to sign.
Note also that the condition that
$\{|\tr \rho(X_n)|\}_n$ is bounded from above
is equivalent to the condition
that the hyperbolic translation lengths of
the isometries $\rho(X_n)$ of $\HH^3$
are bounded from above.

Tan, Wong and Zhang ~\cite{Tan_Wong_Zhang_1, Tan_Wong_Zhang_6}
showed that $\einv(\rho)$
is a closed subset of $\plamination$ and
proved various interesting properties of $\einv(\rho)$,
including a characterization of
those representations $\rho$
with $\einv(\rho)=\emptyset$ or $\plamination$,
generalizing a result of Bowditch ~\cite{Bowditch2}.
They also proposed an interesting conjecture
\cite[Conjecture ~1.8]{Tan_Wong_Zhang_6}
concerning possible homeomorphism types of $\einv(\rho)$.
The following is a modified version of the conjecture
of which Tan ~\cite{Tan} informed the authors.

\begin{conjecture}
\label{conj_TWZ}
{\rm
Suppose $\einv(\rho)$ has at least two accumulation points.
Then either $\einv(\rho)=\plamination$ or
a Cantor set of $\plamination$.
}
\end{conjecture}

They constructed a family of representations $\rho$
which have Cantor sets as $\einv(\rho)$,
and proved the following supporting evidence to the conjecture.

\begin{theorem}
\label{discrete_TWZ}
Let $\rho:\pi_1(\ptorus)\to \SL(2,\CC)$ be {\rm discrete}
in the sense that the set
$\{\tr(\rho(X))\svert X\in\curve\}$
is discrete in $\CC$.
Then if $\einv(\rho)$ has at least three elements,
then $\einv(\rho)$ is either a Cantor set of $\plamination$
or all of $\plamination$.
\end{theorem}

The above theorem implies that
the end invariants $\einv(\rho_r)$
of the representation $\rho_r$ induced by
the holonomy representation of
a hyperbolic $2$-bridge link $K(r)$
is a Cantor set.
But it does not give us the exact description of $\einv(\rho_r)$.
By using the main results stated in Section ~\ref{statements},
we can explicitly determine the end invariants $\einv(\rho_r)$.
To state the theorem,
recall that
the {\it limit set} $\Lambda(\RGPP{r})$
of the group $\RGPP{r}$
is the set of accumulation points in the closure of $\HH^2$
of the $\RGPP{r}$-orbit of a point in $\HH^2$.

\begin{theorem}
For a hyperbolic $2$-bridge link $K(r)$,
the set $\einv(\rho_r)$
is equal to the limit set $\Lambda(\RGPP{r})$
of the group $\RGPP{r}$.
\end{theorem}

We would like to propose the following conjecture.

\begin{conjecture}
{\rm
Let $\rho:\pi_1(\ptorus)\to \PSL(2,\CC)$ be
a type-preserving representation
such that $\einv(\rho)=\Lambda(\RGPP{r})$.
Then $\rho$ is conjugate to the representation $\rho_r$.
}
\end{conjecture}

\section{Application to McShane's identity}
\label{application2}

In his Ph.D. thesis ~\cite{McShane0},
McShane proved the following surprising theorem.

\begin{theorem}
\label{McShane}
Let $\rho:\pi_1(\ptorus)\to \PSL(2,\RR)$
be a type-preserving fuchsian representation.
Then
\[
2\sum_{s\in \QQQ}\frac{1}{1+e^{l_{\rho}(\beta_s)}}=\frac{1}{2}
\]
\end{theorem}

In the above identity, $l_{\rho}(\beta_s)$ denotes the
translation length of
the orientation-preserving isometry $\rho(\beta_s)$
of the hyperbolic plane.
This identity has been generalized
to cusped hyperbolic surfaces by McShane himself ~\cite{McShane},
to hyperbolic surfaces with cusps and geodesic boundary
by Mirzakhani ~\cite{Mirzakhani},
and to hyperbolic surfaces with cusps, geodesic boundary
and conical singularities
by Tan, Wong and Zhang ~\cite{Tan_Wong_Zhang_2}.
A wonderful application to the Weil-Petersson volume of the
moduli spaces of bordered hyperbolic surface was found by
Mirzakhani ~\cite{Mirzakhani}.
Bowditch ~\cite{Bowditch2} (cf. \cite{Bowditch1}) showed
that the identity in Theorem ~\ref{McShane}
is also valid for all quasifuchsian representations of
$\pi_1(\ptorus)$,
where $l_{\rho}(\beta_s)$ is regarded
as the complex translation length of
the orientation-preserving isometry $\rho(\beta_s)$ of the
hyperbolic $3$-space.
Moreover, he gave a nice variation of the identity
for hyperbolic once-punctured torus bundles,
which describes the cusp shape in terms of the
complex translation lengths of essential simple loops
on the fiber torus ~\cite{Bowditch3}.
Other $3$-dimensional variations have been obtained by
\cite{AMS, AMS2, Tan_Wong_Zhang_1, Tan_Wong_Zhang_2, Tan_Wong_Zhang_4,
Tan_Wong_Zhang_5, Tan_Wong_Zhang_6, Tan_Wong_Zhang_7}.

As an application of the main results stated in Section ~\ref{statements},
we can obtain yet another $3$-dimensional variation
of McShane's identity,
which describes the cusp shape of a hyperbolic $2$-bridge link
in terms of the complex translation lengths of
essential simple loops on the bridge sphere.
This proves a conjecture proposed by the first author in \cite{Sakuma}.

To describe the result, note that
each cusp of the hyperbolic manifold $S^3-K(r)$
carries a Euclidean structure, well-defined up to similarity,
and hence it is identified with the quotient of $\CC$
(with the natural Euclidean metric)
by the lattice $\ZZ\oplus \ZZ\lambda$,
generated by the translations
$[\zeta\mapsto \zeta+1]$ and
$[\zeta\mapsto \zeta+\lambda]$
corresponding to the meridian and (suitably chosen) longitude respectively.
This $\lambda$ does not depend on the choice of the cusp,
because
when $K(r)$ is a two-component link
there is an isometry of $S^3-K(r)$ interchanging the two cusps.
We call $\lambda$ the {\it modulus} of the cusp and denote it by
$\lambda(K(r))$.

\begin{theorem}
\label{Conj:McShane}
For a hyperbolic $2$-bridge link $K(r)$ with $r=q/p$,
the following identity holds:
\[
2\sum_{s\in \mathrm{int}I_1}
\frac{1}{1+e^{l_{\rho_r}(\beta_s)}}
+
2\sum_{s\in \mathrm{int}I_2}\frac{1}{1+e^{l_{\rho_r}(\beta_s)}}
+
\sum_{s\in\partial I_1\cup\partial I_2}\frac{1}{1+e^{l_{\rho_r}(\beta_s)}}
=-1.
\]
Further the modulus $\lambda(K(r))$ of the cusp torus of the
cusped hyperbolic manifold $S^3-K(r)$
with respect to a suitable choice of a longitude
is given by the following formula:
\[
\lambda(K(r))=
\begin{cases}
8\sum_{s\in \mathrm{int}I_1}\frac{1}{1+e^{l_{\rho_r}(\beta_s)}}
+
4\sum_{s\in\partial I_1}\frac{1}{1+e^{l_{\rho_r}(\beta_s)}}
&
\mbox{if $p$ is odd,}\\
4\sum_{s\in \mathrm{int}I_1}\frac{1}{1+e^{l_{\rho_r}(\beta_s)}}
+
2\sum_{s\in\partial I_1}\frac{1}{1+e^{l_{\rho_r}(\beta_s)}}
&
\mbox{if $p$ is even.}
\end{cases}
\]
\end{theorem}

The main results stated in Section \ref{statements}
are used to establish the absolute convergence
of the infinite series.

\section{Presentations of 2-bridge link groups}
\label{group_presentation}

\begin{figure}[h]
\begin{center}
\includegraphics{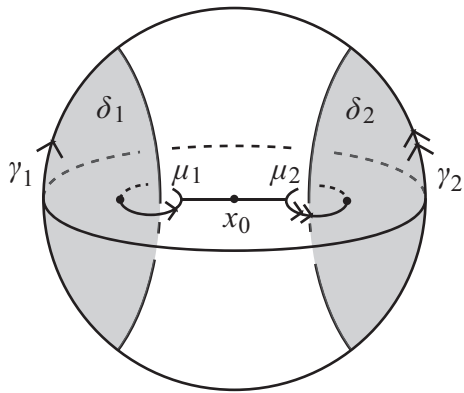}
\end{center}
\caption{
\label{fig.generator}
$\pi_1(B^3-t(\infty), x_0)=F(a,b)$,
where $a$ and $b$ are represented by $\mu_1$ and $\mu_2$, respectively.
}
\end{figure}

In the remainder of this note,
$p$ and $q$ denote relatively prime positive integers
such that $1\le q \le p$ and $r=q/p$.
Theorems ~\ref{first_result} and \ref{second_result} are proved by applying
the small cancellation theory to a
two-generator and one-relator presentation of the link group $G(K(r))$.
To recall the presentation,
let $a$ and $b$, respectively, be the elements of
$\pi_1(B^3-t(\infty), x_0)$
represented by the oriented loops $\mu_1$ and $\mu_2$
based on $x_0$ as illustrated in Figure ~\ref{fig.generator}.
Then $\pi_1(B^3-t(\infty),x_0)$ is identified with
the free group $F(a,b)$.
Note that
$\mu_i$ intersects the disk, $\delta_i$, in $B^3$
bounded by a component of $t(\infty)$ and
the essential arc, $\gamma_i$, on
$\partial(B^3,t(\infty))=\Conways$ of slope $1/0$,
in Figure ~\ref{fig.generator}.

\begin{figure}[h]
\begin{center}
\includegraphics{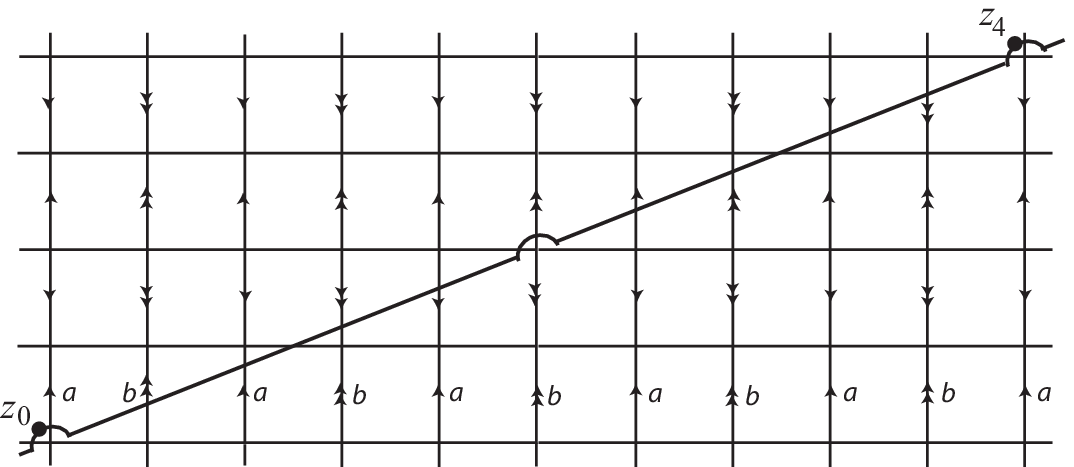}
\end{center}
\caption{\label{line}
The line of slope $2/5$ gives $\hat{u}_{2/5}=bab^{-1}a^{-1}$
and $u_{2/5}=a \hat{u}_{2/5} b \hat{u}_{2/5}^{-1}
=a bab^{-1}a^{-1} b aba^{-1}b^{-1}$.
Since the inverse image of $\gamma_1$ (resp. $\gamma_2$)
in $\RR^2-\ZZ^2$ is the union of
the single arrowed (resp. double arrowed)
vertical edges,
a positive intersection with
a single arrowed (resp. double arrowed) edge
corresponds to $a$ (resp. $b$).
}
\end{figure}

To obtain an element, $u_r$, of $F(a,b)$ represented by the simple loop
$\alpha_r$ (with a suitable choice of an orientation and
a path joining $\alpha_r$
to the base point $x_0$),
note that the inverse image of
$\gamma_1$ (resp. $\gamma_2$)
in $\RR^2-\ZZ^2$ is the union of
the single arrowed (resp. double arrowed)
vertical edges in Figure ~\ref{line}.
Let $L(r)$ be the line in $\RR^2$ of slope $r$
passing through the origin, and let $L^+(r)$ be the line in $\RR^2-\ZZ^2$
obtained by slightly modifying $L(r)$ near each of
the lattice points in $L(r)$
so that $L^+(r)$ takes an upper circuitous route around it,
as illustrated in Figure ~\ref{line}.
Pick a base point $z_0$ from the intersection of $L^+(r)$
with the second quadrant, and consider the sub-line-segment of $L^+(r)$
bounded by $z_0$ and $z_4:=z_0+(2p,2q)$.
Then the image of the sub-line-segment in $\PConway$
is homotopic to the loop $\alpha_s$.
Let $u_r$ be the word in $\{a,b\}$ obtained by reading the
intersection of the line-segment with the vertical lattice lines
(= the inverse images of $\gamma_1$ and $\gamma_2$)
as in Figure ~\ref{line}.
Then $u_r\in F(a,b)\cong\pi_1(B^3-t(\infty))$
is represented by the simple loop $\alpha_r$,
and we obtain the following two-generator one-relator presentation.
\[
\begin{aligned}
G(K(r))&=\pi_1(S^3-K(r))\cong\pi_1(B^3-t(\infty))/\llangle \alpha_r\rrangle \\
&\cong F(a, b)/ \llangle u_r \rrangle
\cong \langle a, b \, | \, u_r \rangle.
\end{aligned}
\]
To describe the explicit formula for $u_r$,
set $\epsilon_i = (-1)^{\lfloor iq/p \rfloor}$
where $\lfloor x \rfloor$ is the greatest integer not exceeding $x$.
Then we have the following (cf. \cite[Proposition ~1]{Riley}).
Let
\[\epsilon_i = (-1)^{\lfloor iq/p \rfloor},\]
where $\lfloor x \rfloor$ is the greatest integer not exceeding $x$.
\begin{enumerate}
\item If $p$ is odd, then
\[u_{q/p}=a\hat{u}_{q/p}b^{(-1)^q}\hat{u}_{q/p}^{-1},\]
where
$\hat{u}_{q/p} = b^{\epsilon_1} a^{\epsilon_2} \cdots b^{\epsilon_{p-2}} a^{\epsilon_{p-1}}$.
\item If $p$ is even, then
\[u_{q/p}=a\hat{u}_{q/p}a^{-1}\hat{u}_{q/p}^{-1},\]
where
$\hat{u}_{q/p} = b^{\epsilon_1} a^{\epsilon_2} \cdots a^{\epsilon_{p-2}} b^{\epsilon_{p-1}}$.
\end{enumerate}
In the above formula, $\hat{u}_{q/p}$ is obtained from the
open interval of $L(r)$ bounded by $(0,0)$ and $(p,q)$.

\begin{figure}[h]
\begin{center}
\includegraphics{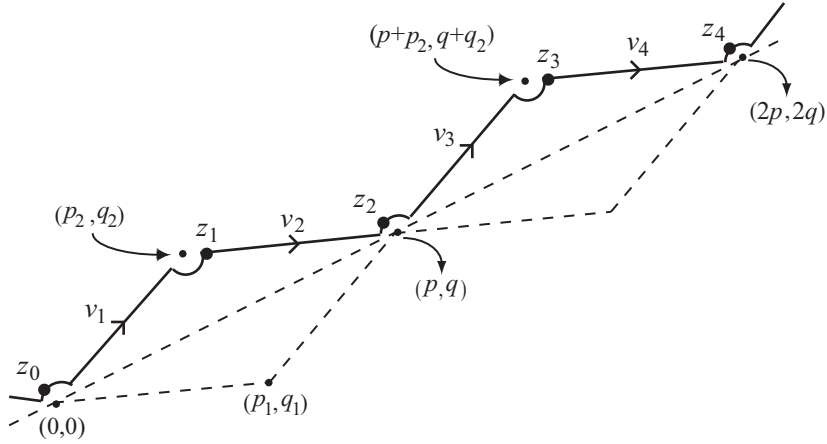}
\end{center}
\caption{\label{line2}
The decomposition of the relator $u_r=v_1v_2v_3v_4$}
\end{figure}

We now describe a natural decomposition of the word $u_r$,
which plays a key role
in the proof of the main results.
Let $r_i=q_i/p_i$ ($i=1,2$) be the rational number
introduced in Section ~\ref{statements}.
Then $(p,q)=(p_1+p_2,q_1+q_2)$
and the parallelogram in $\RR^2$
spanned by $(0,0)$, $(p_1,q_1)$, $(p_2,q_2)$ and $(p,q)$
does not contain lattice points in its interior.
Consider the infinite broken line, $L_b(r)$,
obtained by joining the lattice points
\[
\dots, (0,0), (p_2,q_2), (p,q), (p+p_2, q+q_2), (2p,2q), \dots
\]
which is invariant by the translation $(x,y)\mapsto (x+p,y+q)$.
Let $L_b^+(r)$ be the topological line obtained by
slightly modifying $L_b(r)$ near each of the lattice points in $L_b(r)$
so that $L_b^+(r)$ takes an upper or lower circuitous route around it
according as the lattice point is of the form
$d(p,q)$ or $d(p,q)+(p_2,q_2)$ for some $d\in\ZZ$,
as illustrated in Figure ~\ref{line2}.
We may assume the base points $z_0$ and $z_4$ in $L^+(r)$
also lie in $L_b^+(r)$.
Then the sub-arcs of $L^+(r)$ and $L_b^+(r)$
bounded by $z_0$ and $z_4$ are homotopic in $\RR^2-\ZZ^2$
by a homotopy fixing the end points.
Moreover, the word $u_r$ is also obtained by reading the
intersection of the sub-path of $L_b^+(r)$ with the vertical lattice lines.
Pick a point $z_1\in L_b^+(r)$ whose $x$-coordinate is
$p_2+\mbox{(small positive number)}$, and set
$z_2:=z_0+(p,q)$ and $z_3:=z_1+(p,q)$.
Let $L_{b,i}^+(r)$ be the sub-path of $L_b^+(r)$
bounded by $z_{i-1}$ and $z_i$ ($i=1,2,3,4$),
and consider the subword, $v_i$, of $u_r$ corresponding to $L_{b,i}^+(r)$.
Then we have the decomposition
\[
u_r= v_1v_2v_3v_4,
\]
where the lengths of the subwords $v_i$ are given by
$|v_1|=|v_3|=p_2+1$ and $|v_2|=|v_4|=p_1-1$.
This decomposition plays a key role in the following section.

\section{Sequences associated with the simple loop $\alpha_r$}
\label{sequences}

We begin with the following observation.
\begin{enumerate}
\item
The word $u_r$ is {\it reduced}, i.e.,
it does not contain $xx^{-1}$ or $x^{-1}x$ for any $x \in \{a,b\}$.
It is also {\it cyclically reduced}, i.e., all its cyclic permutations
are reduced.
\item
The word $u_r$ is {\it alternating}, i.e.,
$a^{\pm 1}$ and $b^{\pm 1}$ appear in $u_r$ alternately,
to be precise, neither $a^{\pm2}$ nor $b^{\pm2}$ appears in $u_r$.
It is also  {\it cyclically alternating},
i.e., all its cyclic permutations are alternating.
\end{enumerate}
This observation implies that the word $u_r$ is determined
by the $S$-sequence defined below and the initial letter (with exponent).

\begin{definition}
\label{def:alternating}
{\rm (1) Let $w$ be a nonempty reduced word in
$\{a,b\}$. Decompose $w$ into
\[
w \equiv w_1 w_2 \cdots w_t,
\]
where, for each $i=1, \dots, t-1$, all letters in $w_i$ have positive (resp. negative) exponents,
and all letters in $w_{i+1}$ have negative (resp. positive) exponents.
(Here the symbol $\equiv$ means that the two words are
not only equal as elements of the free group
but also visibly equal, i.e., equal without cancellation.)
Then the sequence of positive integers
$S(w):=(|w_1|, |w_2|, \dots, |w_t|)$ is called the
{\it $S$-sequence of $v$}.

(2) Let $(w)$ be a nonempty reduced cyclic word in
$\{a, b\}$ represented by a word $w$.
Decompose $(w)$ into
\[
(w) \equiv (w_1 w_2 \cdots w_t),
\]
where all letters in $w_i$ have positive (resp. negative) exponents,
and all letters in $w_{i+1}$ have negative (resp. positive) exponents (taking
subindices modulo $t$).
Then the {\it cyclic} sequence of positive integers
$CS(w):=\lp |w_1|, |w_2|, \dots, |w_t| \rp$ is called
the {\it cyclic $S$-sequence of $(w)$}.
Here the double parentheses denote that the sequence is considered modulo
cyclic permutations.
}
\end{definition}
In the above definition,
by a {\it cyclic word},
we mean the set of all cyclic permutations of a
cyclically reduced word.
By $(v)$, we denote the cyclic word associated with a
cyclically reduced word $v$.

\begin{definition}
\label{def4.1(3)}
{\rm
For a rational number $r$ with $0<r\le 1$,
let $u_r$ be the word in $\{a,b\}$ defined in
Section ~\ref{group_presentation}.
Then the symbol $S(r)$ (resp. $CS(r)$) denotes the
$S$-sequence $S(u_r)$ of $u_r$
(resp. cyclic $S$-sequence $CS(u_r)$ of $(u_r)$), which is called
the {\it S-sequence of slope $r$}
(resp. the {\it cyclic S-sequence of slope $r$}).}
\end{definition}

We can easily observe the following.
\[
\begin{aligned}
S(r)
=S(u_r)&=(S(v_1), S(v_2), S(v_3), S(v_4)),\\
CS(r)
=CS(u_r)&=\lp S(v_1), S(v_2), S(v_3), S(v_4)\rp,
\end{aligned}
\]
where $u_r=v_1v_2v_3v_4$ is the natural decomposition of $u_r$
obtained at the end of the last section.
It is also not difficult to observe $S(v_1)=S(v_3)$ and
$S(v_2)=S(v_4)$.
By setting $S_1:=S(v_1)=S(v_3)$ and $S_2:=S(v_2)=S(v_4)$,
we have the following key propositions.

\begin{proposition}
\label{sequence}
The decomposition $S(r)=(S_1, S_2, S_1, S_2)$ satisfies the following.
\begin{enumerate}[\indent \rm (1)]
\item Each $S_i$ is symmetric,
i.e., the sequence obtained from $S_i$ by reversing the order is
equal to $S_i$. (Here, $S_1$ is empty if $k=1$.)

\item Each $S_i$ occurs only twice in
the cyclic sequence $CS(r)$.

\item Set $m:=\lfloor q/p \rfloor$.
Then $S(r)$ consists of only $m$ and $m+1$,
and $S_1$ begins and ends with $m+1$, whereas
$S_2$ begins and ends with $m$.
\end{enumerate}
\end{proposition}

\begin{proposition}
\label{connection}
Let $S(r)= (S_1, S_2, S_1, S_2)$ be as in Proposition ~\ref{sequence}.
For a rational number $s$ with $0 < s \le 1$,
suppose that the cyclic $S$-sequence $CS(s)$ contains
both $S_1$ and $S_2$ as a subsequence.
Then $s \notin I_1 \cup I_2$.
\end{proposition}

\section{Small cancellation conditions for 2-bridge link groups}
\label{small_cancellation}

A subset $R$ of the free group $F(a,b)$ is called {\it symmetrized},
if all elements of $R$ are cyclically reduced and,
for each $w \in R$, all cyclic permutations of $w$ and $w^{-1}$ also belong to $R$.

\begin{definition}
{\rm Suppose that $R$ is
a symmetrized subset of $F(a,b)$.
A nonempty word $v$ is called a {\it piece} if there exist distinct
$w_1, w_2 \in R$ such that $w_1 \equiv vc_1$ and $w_2 \equiv vc_2$.
Small cancellation conditions $C(p)$ and $T(q)$,
where $p$ and $q$ are integers such that $p \ge 2$ and $q \ge 3$,
are defined as follows (see \cite{lyndon_schupp}).
\begin{enumerate}
\item Condition $C(p)$: If $w \in R$
is a product of $n$ pieces, then $n \ge p$.

\item Condition $T(q)$: For
$w_1, \dots, w_n \in R$
with no successive elements $w_i, w_{i+1}$
an inverse pair $(i$ mod $n)$, if $n < q$, then at least one of the products
$w_1 w_2,\dots,$ $w_{n-1} w_n$, $w_n w_1$
is freely reduced without cancellation.
\end{enumerate}
}
\end{definition}

The following proposition enables us to apply small cancellation theory
to the group presentation
$\langle a, b \, |\, u_r \rangle$ of $G(K(r))$.

\begin{proposition}
\label{small_cancellation_condition}
Let $r$ be a rational number such that $0 < r< 1$, and
let $R$ be the symmetrized subset of $F(a, b)$ generated
by the single relator $u_r$
of the group presentation $G(K(r))=\langle a, b \, |\, u_r \rangle$.
Then $R$ satisfies $C(4)$ and $T(4)$.
\end{proposition}

This proposition follows from
the following characterization of pieces,
which in turn is proved by using Proposition ~\ref{sequence}.

\begin{proposition}
\label{cor:max_piece_2}
{\rm (1)} A subword $w$ of the cyclic word
$(u_r^{\pm 1})$ is a piece if and only if
$S(w)$ does not contain $S_1$ as a subsequence
and does not contain $S_2$ in its interior,
i.e., $S(w)$ does not contain a subsequence
$(\ell_1, S_2, \ell_2)$ for some $\ell_1,\ell_2\in\ZZ_+$.

{\rm (2)} For a subword $w$ of the cyclic word $(u_r^{\pm 1})$,
$w$ is not a product of two pieces
if and only if $S(w)$ either contains
$(S_1,S_2)$ as a proper initial subsequence
or contains
$(S_2,S_1)$ as a proper terminal subsequence.
\end{proposition}

\section{Outline of the proof of Theorem ~\ref{first_result}}
\label{outline_proof}

Let $R$ be the symmetrized subset of $F(a, b)$
generated by the single relator $u_r$ of the group presentation
$G(K(r))=\langle a,b \ | \ u_r\rangle$.
Suppose on the contrary that $\alpha_s$ is null-homotopic in
$S^3-K(r)$, i.e., $u_s=1$ in $G(K(r))$, for some $s\in I_1\cup I_2$.
Then there is a {\it van Kampen diagram} $M$ over
$G(K(r))=\langle a, b \, | \, R \, \rangle$
such that the boundary label is $u_s$.
Here $M$ is a simply connected $2$-dimensional complex
embedded in $\RR^2$,
together with a function $\phi$ assigning to
each oriented edge $e$ of $M$, as a {\it label},
a reduced word $\phi(e)$ in $\{a,b\}$ such that the following hold.
\begin{enumerate}
\item If $e$ is an oriented edge of $M$ and $e^{-1}$ is the oppositely oriented edge,
then $\phi(e^{-1})=\phi(e)^{-1}$.
\item For any boundary cycle $\delta$ of any face of $M$,
$\phi(\delta)$ is a cyclically reduced word representing
an element of $R$.
(If $\alpha=e_1, \dots, e_n$ is a path in $M$, we define $\phi(\alpha) \equiv \phi(e_1) \cdots \phi(e_n)$.)
\end{enumerate}

We may assume $M$ is {\it reduced},
namely it satisfies the following condition:
Let $D_1$ and $D_2$ be faces (not necessarily distinct) of $M$
with an edge $e \subseteq \partial D_1 \cap \partial D_2$,
and let $e \delta_1$ and $\delta_2e^{-1}$ be boundary cycles of $D_1$ and $D_2$, respectively.
Set $\phi(\delta_1)=f_1$ and $\phi(\delta_2)=f_2$.
Then we have $f_2\ne f_1^{-1}$.
Moreover, we may assume the following conditions:
\begin{enumerate}
\item
$d_M(v) \ge 3$ for every vertex $v \in M-\partial M$.
\item
For every edge $e$ of $\partial M$, the label $\phi(e)$ is a piece.
\item
For a path $e_1, \dots, e_n$ in $\partial M$ of length $n\ge 2$
such that the vertex
$e_i\cap e_{i+1}$ has degree $2$
for $i=1,2,\dots, n-1$,
$\phi(e_1) \phi(e_2) \cdots\phi(e_n)$ cannot be expressed as a product of less than $n$ pieces.
\end{enumerate}

Since $R$ satisfies the conditions $C(4)$ and $T(4)$
by Proposition ~\ref{small_cancellation_condition},
$M$ is a $[4,4]$-map, i.e.,
\begin{enumerate}
\item $d_M(v) \ge 4$ for every vertex $v \in M-\partial M$;
\item $d_M(D) \ge 4$ for every face $D \in M$.
\end{enumerate}
Here, $d_M(v)$, the {\it degree of $v$},
denotes the number of oriented edges in $M$ having $v$ as initial vertex,
and $d_M(D)$, the {\it degree of $D$}, denotes the number of
oriented edges in a boundary cycle of $D$.

Now, for simplicity, assume that $M$ is homeomorphic to a disk.
(In general, we need to consider an extremal disk of $M$.)
Then by the Curvature Formula of
Lyndon and Schupp (see ~\cite[Corollary ~V.3.4]{lyndon_schupp}),
we have
\[
\sum_{v \in \, \partial M} (3-d_M(v)) \ge 4.
\]
By using this formula, we see that
there are three edges $e_1$, $e_2$ and $e_3$ in $\partial M$
such that $e_1 \cap e_2=\{v_1\}$
and $e_2 \cap e_3=\{v_2\}$, where $d_M(v_i)=2$ for each $i=1, 2$.
Since $\phi(e_1)\phi(e_2)\phi(e_3)$ is not expressed as a
product of two pieces, we see
by Proposition ~\ref{cor:max_piece_2}
that the boundary label of $M$
contains a subword, $w$, with $S(w)=(S_1,S_2,\ell)$ or $(\ell, S_2,S_1)$.
This in turn implies that
the $S$-sequence of the boundary label contains both $S_1$ and $S_2$
as subsequences.
Hence, by Proposition ~\ref{connection},
we have $s\not\in I_1\cup I_2$, a contradiction.

\section{Outline of the proof of Theorem ~\ref{second_result}}
\label{outline_proof2}

Suppose, for two distinct $s, s' \in I_1\cup I_2$,
the unoriented loops $\alpha_s$ and $\alpha_{s'}$ are homotopic in
$S^3-K(r)$.
Then there is a reduced annular $R$-diagram $M$
such that $u_s$ is an outer boundary label and
$u_{s'}^{\pm 1}$ is an inner boundary label of $M$.
Again we can see that $M$ is a $[4,4]$-map
and hence we have the following curvature formula:
\[
\sum_{v \in \partial M} (3-d_M(v)) \ge 0.
\]
By using this formula,
we obtain the following very strong structure theorem
for $M$, which plays key roles throughout the series
of papers \cite{lee_sakuma_2, lee_sakuma_3, lee_sakuma_4}.

\begin{figure}[h]
\begin{center}
\includegraphics{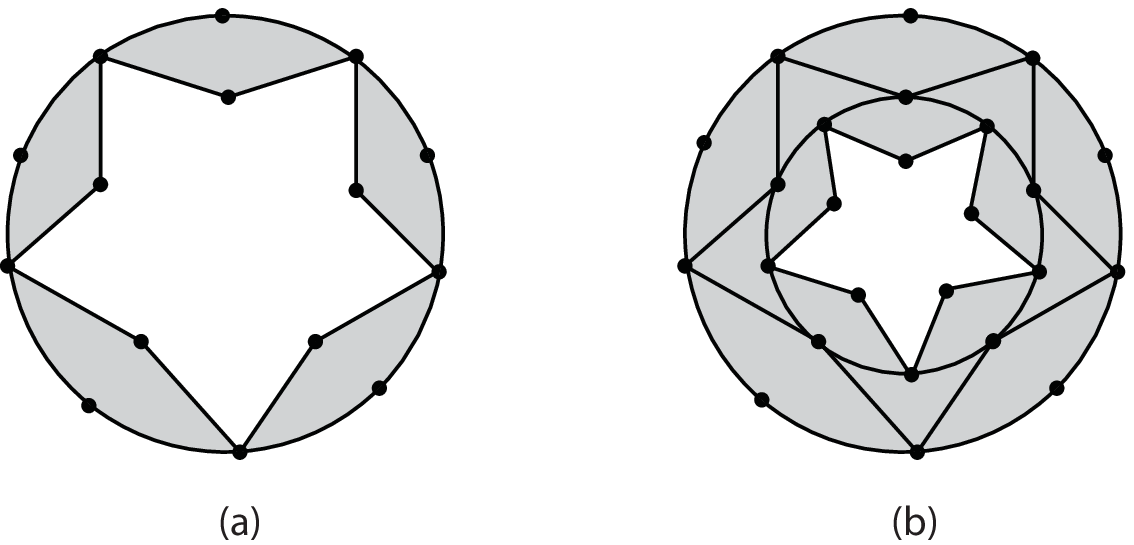}
\end{center}
\caption{\label{layer}}
\end{figure}

\begin{theorem}
\label{cor:structure}
Figure ~\ref{layer}(a) illustrates the only possible type of
the outer boundary layer of $M$,
while Figure ~\ref{layer}(b) illustrates
the only possible type of whole $M$.
(The number of faces per layer and the number of layers are variable.)
\end{theorem}

In the above theorem,
the {\it outer boundary layer} of
the annular map $M$ is the submap of $M$
consisting of all faces $D$
such that the intersection of $\partial D$ with the outer boundary of $M$ contains an edge,
together with the edges and vertices contained in $\partial D$.

The first paper ~\cite{lee_sakuma_2}
of the series devoted to proof of Theorem ~\ref{second_result}
treats the case
when the 2-bridge link is a $(2,p)$-torus link,
the second paper ~\cite{lee_sakuma_3} treats the case
of 2-bridge links of slope $n/(2n+1)$ and $(n+1)/(3n+2)$,
where $n \ge 2$ is an arbitrary integer,
and the third paper ~\cite{lee_sakuma_4} treats the general case.
The two families treated in the second paper
play special roles
in the project in the sense that
the treatment of these links form a base step
of an inductive proof of the theorem for general $2$-bridge links.
We note that the figure-eight knot is both
a 2-bridge link of slope $n/(2n+1)$ with $n=2$ and a 2-bridge link of slope $(n+1)/(3n+2)$ with $n=1$.
Surprisingly,
the treatment of the figure-eight knot,
the simplest hyperbolic knot, is
the most complicated.
This reminds us of the phenomenon
in the theory of exceptional
Dehn filling that
the figure-eight knot attains the
maximal number of exceptional Dehn fillings.

\bibstyle{plain}

\end{document}